\documentclass[titlepage,11pt]{article}
\usepackage{latexsym}
\usepackage{amsfonts}
\usepackage{amsmath}
\newcommand\blackslug{\hbox{\hskip 1pt \vrule width 4pt height 8pt depth 1.5pt
        \hskip 1pt}}
\newcommand\bbox{\hfill \quad \blackslug \bigbreak}
\def\DD{\hbox{-}}
\def\CC{\hbox{-}\cdots\hbox{-}}

\def\cupcup{\cup\cdots\cup}
\def\idw{\operatorname{idw}}
\def\odw{\operatorname{odw}}
\def\glc{\operatorname{glc}}
\def\mcw{\operatorname{mcw}}
\def\ad{\operatorname{ad}}

\title{Bounded-diameter tree-decompositions}
\author{Eli Berger\\
University of Haifa, Haifa, Israel
\\
\\
Paul Seymour\thanks{Supported by AFOSR grant
FA9550-22-1-0234, and NSF grant  DMS-2154169.}\\
Princeton University, Princeton, NJ 08544}

\date{January 4, 2023; revised \today}

\newtheorem{thm}{}[section]

\newcommand{\Proof}{\noindent{\bf Proof.}\ \ }

\begin{document}
\maketitle
\begin{abstract}
When does a graph admit a tree-decomposition in which every bag has small diameter?
For finite graphs, this is a property of interest in algorithmic graph theory, where it is called having bounded ``tree-length''.
We will show that this is equivalent to 
being ``boundedly quasi-isometric to a tree'', which for infinite graphs is a much-studied property from metric geometry. 
One object of this paper is to tie these two areas together. We will prove that
there is a tree-decomposition in which each bag has small diameter, if and only if
there is a map $\phi$ from $V(G)$ into the vertex set of a tree $T$,
such that for all $u,v\in V(G)$, the distances $d_G(u,v), d_T(\phi(u),\phi(v))$
differ by at most a constant.

A necessary condition for admitting such a tree-decomposition is that there is no long geodesic cycle, and for graphs of 
bounded tree-width,
Diestel and M\"uller showed that this is also sufficient. But it is not sufficient in general, even qualitatively, 
because there are graphs
in which every geodesic cycle has length at most three, and yet every tree-decomposition has a bag with large diameter. 

There is a more general necessary condition, however. A ``geodesic loaded cycle'' 
in $G$ is a pair $(C,F)$, where $C$ is a cycle of $G$ and 
$F\subseteq E(C)$, such that for every pair 
$u,v$ of vertices of $C$, 
one of the paths
of $C$ between $u,v$ contains at most $d_G(u,v)$ $F$-edges, where $d_G(u,v)$ is the distance between $u,v$ in $G$.
We will show that a (possibly infinite) graph $G$ admits a tree-decomposition in which every bag has small diameter, if and only if $|F|$ is small for every 
geodesic loaded cycle $(C,F)$. Our proof is an extension of an algorithm to approximate tree-length in finite graphs by Dourisboure and Gavoille.

In metric geometry, there is a similar theorem that characterizes when a graph is quasi-isometric to a tree, 
``Manning's bottleneck criterion''. 
The goal of this paper is to tie all these concepts together, and add a few more related ideas. For instance,
we prove a conjecture of Rose McCarty, that $G$
admits  a tree-decomposition in which every bag has small diameter, if and only if for all vertices $u,v,w$ of $G$, some ball of small radius meets every path
joining two of $u,v,w$. 

\end{abstract}

\section{Introduction}
Graphs in this paper may be infinite.
(Our research was motivated by interest in finite graphs, but all the proofs work equally well for infinite graphs.)
A {\em tree-decomposition} of a graph $G$ is a pair $(T,(B_t:t\in V(T)))$, where $T$ is a tree, and $B_t$
is a subset of $V(G)$ for each $t\in V(T)$, such that:
\begin{itemize}
\item $V(G)$ is the union of the sets $B_t\;(t\in V(T))$; 
\item for every edge $e=uv$ of $G$, there exists $t\in V(T)$ with $u,v\in B_t$; and
\item for all $t_1,t_2,t_3\in V(T)$, if $t_2$ lies on the path of $T$ between $t_1,t_3$, then 
$B_{t_1}\cap B_{t_3}\subseteq B_{t_2}$.
\end{itemize}
($T$ might be infinite.)
The {\em width} of a tree-decomposition $(T,(B_t:t\in V(T)))$ is the maximum of the numbers $|B_t|-1$ for $t\in V(T)$,
or $\infty$ if there is no finite maximum;
and the {\em tree-width} of $G$ is the minimum width of a tree-decomposition of $G$.

We have a good grasp of what stops a graph having bounded tree-width: 
\begin{thm}\label{RS5}
{\bf Theorem {\cite{RS5}}:} There is a function $f$ such that for every graph $G$, if $k\ge 2$ is maximum such that $G$ contains the $k\times k$ grid
as a minor, then the tree-width of $G$ is between $k$ and $f(k)$.
\end{thm}
Indeed, in~\cite{chuzhoy, chuzhoy2} it is shown that $f$ can be chosen to be a polynomial. 

But what if we want a tree-decomposition $(T,(B_t:t\in V(T)))$ such that $G[B_t]$ is connected for each $t\in V(T)$?
Then the requisite size of the bags $B_t$ (the {\em connected tree-width}) may change dramatically. For instance, if $G$ is a cycle of length $\ell$,
then its tree-width is two, but if we want all bags to induce connected subgraphs, then some bags must have
size at least $\ell/3+1$. (This follows from \ref{cycledecomp}, taking $F=E(G)$.) Let us say a connected subgraph $C$ of $G$ is {\em geodesic} if for every two vertices $u,v\in V(C)$,
the distance between $u,v$ in $G$ equals their distance in $C$. So, if $G$ has a geodesic cycle $C$, then its connected tree-width
is at least $|C|/3$, and at least the tree-width. Diestel and M\"uller showed a beautiful converse:
\begin{thm}\label{diestel} 
{\bf Theorem \cite{diestel}:} There is a function $f$ such that if a graph $G$ has tree-width at most $w$ and its longest geodesic cycle has length $\ell$
then its connected tree-width is at most $f(w,\ell)$.
\end{thm}

What happens to \ref{diestel} if we drop the assumption of bounded tree-width, and just assume there is no long geodesic cycle?
Is this qualitatively equivalent to some sort of decomposition? ``Admitting a tree-decomposition in which every 
bag is a connected subgraph of bounded size'' is too strong (for instance, because $G$ might be a large complete graph), 
and ``admitting a tree-decomposition in which every   
bag is a connected subgraph'' is too weak (because every connected graph has such a tree-decomposition, with a one-vertex tree).
What about ``admitting a tree-decomposition in which every
bag is a connected subgraph of bounded diameter''? In one direction this works: if $G$ has a geodesic cycle of length at least 
$\ell$,  then 
in every tree-decomposition, some bag has diameter at least $\ell/3$.
What about the converse? Is it true that if $G$ 
does {\em not} have a long geodesic cycle, then it admits a tree-decomposition such that each bag has bounded diameter?
At first sight this looks plausible. For instance, geodesic cycles are induced, and if we replace ``no long geodesic cycle'' with ``no long
induced cycle'' the result is true for finite graphs (and presumably also for infinite ones, though we have not checked).
\begin{thm}\label{uncle}
{\bf Theorem \cite{uncle}:} For all integers $\ell\ge 4$, if $G$ is a finite graph with no induced cycle of length $>\ell$,
then $G$ admits a tree-decomposition $(T,(B_t:t\in V(T)))$ such that for each $t\in V(T)$, every two vertices of $B_t$
are joined by a path of $G[B_t]$ with length at most $\ell$.
\end{thm}

But the attempt at a converse proposed above is wrong: as we shall see in \ref{geodesic}, there are graphs with no geodesic cycle of length more than three, in which every tree-decomposition
has a bag of large diameter. So now there are two questions: what is the right ``structural'' statement that goes with
not having a long geodesic cycle; and what is the right ``exclusion'' statement that goes with admitting a tree-decomposition 
with bags of small diameter? We have not answered the first question, but we can answer the second, and that is the primary
goal of this paper.

It turns out that
the right thing to exclude is indeed a kind of cycle, a ``geodesic loaded cycle''.
Let us be more precise. 
If $u,v\in V(G)$, then $d_{G}(u,v)$ denotes the length (that is, number of edges) of the shortest path of $G$
between $u,v$, or $\infty$ if there is no such path.
Let $C$ be a cycle of $G$ and let $F\subseteq E(C)$. We call the pair $(C,F)$ a {\em loaded cycle} of $G$, and $|F|$ is its {\em load}.
If $u,v\in V(C)$ are distinct, let $d_{C,F}(u,v)$ denote the smaller of 
$|E(P)\cap F|,|E(Q)\cap F|$ where $P,Q$ are the two paths of $C$ between $u,v$. (Let $d_{C,F}(u,v):=0$ if $u=v$.)
Let us say that the loaded cycle $(C,F)$ is {\em geodesic} in $G$ if 
$d_G(u,v)\ge d_{C,F}(u,v)$ for all $u,v\in V(C)$. If $G$ admits a tree-decomposition in which all bags have bounded diameter,
then every geodesic loaded cycle has bounded load; and our main theorem says that if every geodesic loaded cycle has bounded 
load, then $G$ admits a tree-decomposition in which all bags have bounded diameter.

Incidentally, what does ``all bags have bounded diameter'' mean? We might mean that for each bag, every two of its vertices
are at bounded distance in $G$; or we might mean that for each bag, every two of its vertices are at bounded distance in
the subgraph induced on the bag. Fortunately, these two turn out to be essentially equivalent, in the sense that if $G$ 
admits a tree-decomposition
in which each bag has diameter at most $d$ (measuring distance in $G$), then $G$ also admits a tree-decomposition
in which each bag has diameter at most $2d$ (measuring distance in the bag).

We need some more definitions.
Let $\glc(G)$ be the maximum load over all geodesic loaded cycles in $G$, or $\infty$ if there is no such maximum,
or $0$
if $G$ has no geodesic loaded cycle (and hence $G$ has no cycle).
If $\mathcal{B}=(T,(B_t:t\in V(T)))$ is a tree-decomposition of $G$, we define
the {\em inner diameter}  of $\mathcal{B}$ to be the maximum of the diameter of $G[B_t]$ for $t\in V(T)$ (and so $\infty$ if 
some $G[B_t]$ is not connected, or its diameter is unbounded);
and we define the {\em inner diameter-width} $\idw(G)$ of $G$ to be the minimum of the inner diameter  over all tree-decompositions 
$G$.
Similarly, 
we define the {\em outer diameter} of $\mathcal{B}$ to be
$$\max_{t\in V(T)} \max_{u,v\in B_t} d_{G}(u,v),$$
(if it exists, and $\infty$ otherwise) and the {\em outer diameter-width} $\odw(G)$ of $G$ to be the minimum of the outer diameter of $\mathcal{B}$ over all tree-decompositions $\mathcal{B}$ of $G$.
(Note that we are not bothering with the customary $-1$ in these definitions of width.)
Outer diameter-width is called ``tree-length'' in algorithmic graph theory~\cite{dissaux,dourisboure}, but we stick with ``outer diameter-width'' here to 
emphasize the distinction with ``inner diameter-width''.
We will show that these three numbers are related, with the following two theorems:

\begin{thm}\label{onetotwo0}
{\bf Theorem:} For every graph $G$, $\odw(G)\le \idw(G)\le 2\odw(G)$.
\end{thm}
\begin{thm}\label{mainthm}
{\bf Theorem:} For every graph $G$, $\odw(G)-1\le \glc(G)\le 3\odw(G)$.
\end{thm}
In section~\ref{sec:conseq} we will prove \ref{onetotwo0} and the second inequality of \ref{mainthm}, and in section \ref{sec:mainthm} we will 
prove the rest of \ref{mainthm}.
None of these is difficult, but the first inequality of \ref{mainthm} is the least easy. It is closely related to a theorem
of Manning~\cite{manning} in metric geometry, and to an approximation algorithm of Dourisboure and Gavoille~\cite{dourisboure}, as we will explain
later.

Let $T$ be a tree, and let $\phi$ be a map from $V(G)$ into $V(T)$. The {\em additive distortion} of 
$(T,\phi)$ is the maximum of  
$$|d_G(u,v)-d_T(\phi(u),\phi(v))|$$
over all $u,v\in V(G)$ (or $\infty$ if this is unbounded).
The {\em additive distortion} $\ad(G)$ of $G$ is the minimum 
$k$ such that there is a tree $T$ and a map $\phi:V(G)\to V(T)$ with additive distortion at most $k$.
We will prove:
\begin{thm}\label{kerrapp}
{\bf Theorem:} Let $G$ be a connected graph. Then $(\odw(G)-1)/2\le \ad(G)\le 6\odw(G)+1$.
\end{thm}
The connection between outer diameter-width and additive distortion seems to be new, and exploring it is
a second goal of this paper.

There is a more general relation, ``quasi-isometry''. 
(This is a concept from metric spaces, but we will define it just for graphs.)
Let $G,H$ be graphs, and let $\phi:V(G)\to V(H)$ be a map.
Let $L\ge 1$ and $C\ge 0$; we say that $\phi$ is an {\em $(L,C)$-quasi-isometry}
if:
\begin{itemize}
\item for all $u,v$ in $V(G)$,
$\frac1L d_G(u,v)-C\le d_H(\phi(u),\phi(v))\le L d_G(u,v)+C$;
and
\item for every $y\in V(H)$ there exists $v\in V(G)$ such that $d_H(\phi(v), y)\le  C$.
\end{itemize}
For a connected graph $G$, there is a $(1,C)$-quasi-isometry to a tree
if and only if $\ad(G)\le C$; so 
quasi-isometry to a tree looks more general than additive distortion, because $L$ might be bigger than 1.
But it is not really more general, since a theorem of Kerr implies:
\begin{thm}\label{kerr}
{\bf Theorem \cite{kerr}:} For all $L,C$ there exists $C'$ such that if there is an $(L,C)$-quasi-isometry from a graph $G$ to a tree, then there
is a $(1,C')$-quasi-isometry from $G$ to a tree.
\end{thm}

In section \ref{sec:quasi}, we will prove a result that contains Kerr's theorem (for graphs; Kerr's theorem is really about metric spaces):
\begin{thm}\label{quasitobags0}
{\bf Theorem:} If there is an $(L,C)$-quasi-isometry from a graph $G$ to a tree,  then $G$ is connected and $\odw(G)\le L(L+C+1)+C$.
Conversely, for every connected graph $G$ with $\odw(G)$ finite, there is a $(1,6\odw(G))$-quasi-isometry to a tree.
\end{thm}

There are more graph parameters that are related to outer diameter-width.
Let us say a graph $G$ has {\em McCarty-width} $k$ if $k\ge 0$ is minimum such that the following holds:
for every three vertices $u,v,w$ of $G$, there is a vertex $x$, such that if $X$ denotes the set of all vertices that have distance
at most $k$ from $x$, then no component of $G\setminus X$ contains two of $u,v,w$. 
Let $\mcw(G)$ denote the McCarty-width of $G$; and $\mcw(G)=\infty$ if there is no such $k$.
Rose McCarty suggested that $\odw(G)$ is small
if and only if $\mcw(G)$ is small. This turns out to be true, because of the following, which we will prove in
section~\ref{sec:mccarty}:

\begin{thm}\label{mccarty0}
{\bf Theorem:} Let $G$ be a graph. Then $(\odw(G)-3)/6\le \mcw(G)\le \odw(G)$.
\end{thm}

Finally, in section \ref{sec:counterex} we discuss a different (but false) 
hope for a characterization of when $\odw(C)$ is bounded. If a connected graph $G$ admits a tree-decomposition with small inner 
diameter,
one might hope to ``approximate'' $G$ by a spanning tree  --- is there necessarily a spanning tree $T$ such that all 
distances in $T$ are about the same as the corresponding distance in $G$? The answer is no; there are finite graphs $G$
with $\idw(G)=1$, such that for every spanning tree $T$, there is an edge $uv$ of $G$ with $d_T(u,v)$ arbitrarily large.

\section{Consequences of bounded outer diameter-width}\label{sec:conseq}

We need the following basic fact about tree-decompositions:
\begin{thm}\label{intercept}
{\bf Lemma: } If $(T,(B_t:t\in V(T)))$ is a tree-decomposition of $G$, and $r,s,t\in V(T)$, and $s$ lies in the path of $T$ between $r,t$,
then every path of $G$ with one end in $B_r$ and the other in $B_t$ has a vertex in $B_s$.
\end{thm}
\Proof
If $F$ is a non-null connected subgraph of $G$, then $\{t\in V(T): V(F)\cap B_t\ne \emptyset\}$ is the vertex set of a subtree of $T$ (this
is easily proved by induction on $V(F)$); and the result follows by letting $F$ be a path between $B_r$ and $B_t$.~\bbox

We begin with \ref{onetotwo0}, which we restate:
\begin{thm}\label{onetotwo}
{\bf Theorem:} For every graph $G$, $\odw(G)\le \idw(G)\le 2\odw(G)$.
\end{thm}
\Proof
Clearly $\odw(G)\le \idw(G)$, and we need to prove the second inequality.
Let $(T,(B_t:t\in V(T)))$ be a tree-decomposition of $G$ with outer diameter $\odw(G)$, and let $d:=\odw(G)$. 
We may assume that $d$ is finite.
If $X\subseteq V(G)$, let us define $X^+$ to be the union of the vertex sets of all paths $P$ of $G$ with length at most $d$
and with ends in $X$ (thus $X\subseteq X^+$, because $P$ is permitted to have only one vertex). We claim that
$(T,(B_t^+:t\in V(T)))$ is a tree-decomposition of $G$, and to show this, we only need to show the following:
\\
\\
{\bf (1) Claim: }{\em If $r,s,t\in V(T)$ and $s$ belongs to the path of $T$ between $r,t$, and $w\in V(G)$ belongs to both $B_r^+$
and $B_t^+$, then $w\in B_s^+$.}
\\
\\
Choose $q\in V(T)$ such that $w\in B_q$. Since $s$ belongs to the path of $T$ between $r,t$, it follows that either
$s$ belongs to the path of $T$ between $q,r$, or $s$ belongs to the path of $T$ between $q,t$, and without loss of generality we
may assume the latter. Since $w\in B_t^+$, there is a path $P$ of $G$ with length at most $d$ and with $w\in V(P)$,
such that the ends of $P$ belong to $B_t$. Let the ends of $P$ be $p_1,p_2$ (possibly $p_1=p_2$, if $P$ has length zero), and for $i = 1,2$ let $Q_i$
be the subpath of $P$ between $w,p_i$. Since $V(Q_i)\cap B_t\ne \emptyset$ and $V(Q_i)\cap B_q\ne \emptyset$, and
$s$ belongs to the path of $T$ between $q,t$, it follows that there exists $s_i\in V(Q_i)\cap B_s$ by \ref{intercept},
for $i = 1,2$. Then the subpath of $P$ between $s_1,s_2$ has length at most $d$, and has both ends in $B_s$, and contains
$w$, and so $w\in B_s^+$. This proves (1).

\bigskip

We claim that the tree-decomposition $(T,(B_t^+:t\in V(T)))$ has inner diameter-width at most $2d$. To see this, let
$t\in V(T)$, and let $u,v\in B_t^+$. There is a path $P$ of $G$ with length at most $d$ and with ends in $B_t$
that contains $u$; and so $V(P)\subseteq B_t^+$, and there is a subpath $P'$ of $P$ between $u$ and some vertex $u'\in B_t$
that has length at most $d/2$. Similarly there is a path $Q'$ with $V(Q')\subseteq B_t^+$, of length at most $d/2$,
between $v$ and some vertex  $v'\in B_t$. But there is a path $R$ of $G$ between $u', v'$ of length at most $d$, since
$(T,(B_t:t\in V(T)))$ has outer diameter-width $d$; and so $V(R)\subseteq B_t^+$. The union of $P',Q'$ and $R$ contains
a path between $u,v$ of length at most $2d$ with all vertices in $B_t^+$. This proves that
$(T,(B_t^+:t\in V(T)))$ has inner diameter-width at most $2d$, and so $\idw(G)\le 2\odw(G)$. This proves \ref{onetotwo}.~\bbox

Let us turn to \ref{mainthm}. We need the following lemma:

\begin{thm}\label{cycledecomp}
{\bf Lemma: }Let $C$ be a cycle of a graph $G$, let $F\subseteq E(C)$ with $|F|\ge 2$, and let $(T,(B_t:t\in V(T)))$ be a tree-decomposition of $G$. Then there
exist $t\in V(T)$ and $u,v\in V(C)\cap B_t$ such that $d_{C,F}(u,v)\ge |F|/3$.
\end{thm}
\Proof Suppose not. Since $|B_t\cap V(C)|\ge 2$ for some $t\in V(T)$, it follows that $|F|\ge 4$. If $P$ is a path or cycle, let us say its {\em $F$-length}
is $|F\cap E(P)|$.
\\
\\
{\bf (1) Claim: }{\em For each $t\in V(T)$, there is a path $P_t$ of $C$ with $F$-length less than $|F|/3$, such that $V(C)\cap B_t\subseteq V(P_t)$.}
\\
\\
We may assume that $|V(C)\cap B_t| \ge 2$. By choosing distinct $u,v\in V(C)\cap B_t$ and the supposed falsity of the theorem,
we deduce that there is a path $P$ of $C$, with distinct ends both in $B_t$, and with $F$-length less than $|F|/3$.
Choose such a path $P$ such that $B_t\cap V(P)$ is maximal. Let $P$ have ends $u,v$ say. Suppose that there
exists $w\in V(C)\cap B_t$ with $w\notin V(P)$. By the falsity of the theorem, there are paths $Q,R$ of $C$,
both with $F$-length less than $|F|/3$, joining $u,w$ and $v,w$ respectively. Since $P,Q,R$ all contain fewer than $|F|/3$ edges in $F$,
there is an edge $e$ of $F$ that belongs to none of $P,Q,R$. But then $P,Q,R$ are subpaths of the path
$C\setminus \{e\}$, and so one of $Q,R$ includes $P$, contrary to the maximality of $P$. Thus there is
no such $w$. This proves (1).

\bigskip
For every subtree $T'$ of $T$, let $B(T')$ denote $\bigcup_{s\in V(T')}B_s$.
For each $t\in V(T)$, since $|F|\ge 2$, and more than $2|F|/3$ edges in $F$ do not belong to $P_t$,
there are at least two edges of $F$ that do not belong to $P_t$;
and in particular, $C\setminus V(P_t)$ is a path $Q_t$ say. Since $B_t\cap V(Q_t)=\emptyset$, and $Q_t$ is non-null 
and connected, \ref{intercept} implies that 
there is a component $T_t$ of $T\setminus \{t\}$ such that $V(Q_t)\subseteq B(T_t)$.

Let $P$ be a path of $C$, maximal such that $P=P_t$ for some $t\in V(T)$. 
Let $w\in V(C)\setminus V(P)$, let $r\in V(T)$ with $w\in B_r$, and choose $t\in V(T)$
with $P_t=P$ such that $d_T(r,t)$ is as small as possible. Let $t'$ be the neighbour of $t$ in $T_t$.
Let $P$ have ends $u,v$, and let
$uu'$ be an edge of
$C$ that does not belong to $E(P)$. Choose $s\in V(T)$ with $u,u'\in B_s$. 
Since $u'\in V(Q_t)\cap B_s$, and $u'\notin B_t$, it follows from \ref{intercept} that
$s\in V(T_t)$. Since $u\in B_s$ and $u\in B_t$, we deduce that $u\in B_{t'}$,
and similarly $v\in B_{t'}$. Consequently $P_{t'}$ includes a path of $C$ between $u,v$,
and so includes $P$ (it cannot include the other path of $C$ between $u,v$ since that contains more than $2|F|/3$
edges in $F$). From the maximality of $P$, it follows that $P_{t'}=P$. But $w\in V(Q_t)$ and $w\notin B_t$, so $r\in V(T_t)$
by \ref{intercept}; and consequently  $d_T(r,t')< d_T(r,t)$, contrary to the choice of $t$.
This proves \ref{cycledecomp}.~\bbox

We deduce half of \ref{mainthm}, the following:

\begin{thm}\label{mainthm1}
{\bf Lemma: }For every graph $G$, $\glc(G)\le 3\odw(G)$.
\end{thm}
\Proof We may assume that $\glc(G)\ge 1$; so $G$ has an edge, and so $\odw(G)\ge 1$, and hence we may assume that $\glc(G)>3$.
We may also assume that $\odw(G)$ is finite.
Let $(T,(B_t:t\in V(T)))$ be a tree-decomposition of $G$ with outer diameter-width equal to $\odw(G)$. Let $(C,F)$ be a geodesic
loaded cycle of $G$ with $|F|\ge 2$.
By \ref{cycledecomp}, there
exist $t\in V(T)$ and $u,v\in V(C)\cap B_t$ such that $d_{C,F}(u,v)\ge |F|/3$. Since $(C,F)$ is geodesic, it follows that
$d_G(u,v)\ge |F|/3$; but $d_G(u,v)\le \odw(G)$, and so $\odw(G)\ge |F|/3$. This proves \ref{mainthm1}.~\bbox

\ref{cycledecomp} has another useful consequence. We mentioned earlier that having no long geodesic cycle was not sufficient for $\odw(G)$ to be small; let us prove that.

\begin{thm}\label{geodesic}
{\bf Theorem: }There are finite graphs $G$ with $\odw(G)$ arbitrarily large, in which every geodesic cycle has length three.
\end{thm}
\Proof
Take a large triangular piece of the triangular lattice. Thus, 
let $n$ be some large number, and let $V$ be the set of all triples $(a,b,c)$ of nonnegative integers such that $a+b+c=n$.
We make $(a,b,c)$ and $(a',b',c')$ adjacent if $|a-a'|+|b-b'|+|c-c'|=2$. Let $G$ be the graph just made, and suppose that $C$
is a geodesic cycle of $G$. For $0\le i\le n$, let $P_i$ be the path 
$$(0,n-i,i)\DD(1,n-i-1,i)\CC (n-i,0,i).$$
Then $P_i$ is 
a geodesic path, and for any two vertices in it, the subpath of $P_i$ between them is the only geodesic of $G$ between them.
Consequently, if there are two vertices of $P_i$ in $C$, then $C$ also contains the subpath of $P_i$ between them;
and so the intersection of $P_i$ with $C$ is connected. It follows that $C\setminus V(P_i\cap C)$ is also connected, and lives
completely on one side of $P_i$ in the drawing. The same is true for the paths 
$$(i, 0, n-i)\DD(i,1,n-i-1)\CC (i,n-i,0)$$
and 
$$(0,i,n-i)\DD (1,i,n-i-1)\CC (n-i,i,0).$$ 
Call these $Q_i$ and $R_i$ respectively. Since $C$ is not separated by any of the paths
$P_i,Q_i, R_i$ for $0\le i\le n$, it follows that $V(C)$ belongs to a region of the drawing formed by the union of these paths,
and so has length three. 

Next we need to show that $\odw(G)$ is large. Let $C$ be the 
perimeter cycle of $G$, of length $3n$. 
It is easy to check (and we omit the details) that $(C,E(P_0))$ is a geodesic loaded cycle,
so by \ref{mainthm1}, 
$$\odw(G)\ge \glc(G)/3\ge |E(P_0)|/3=n/3.$$ 
This proves \ref{geodesic}.~\bbox 
\section{Outer diameter-width and geodesic loaded cycles} \label{sec:mainthm}

If $P$ is a path, its {\em interior} $P^*$ is the set of vertices of $P$ that have degree two in $P$.
Now we prove the remainder of \ref{mainthm}, in the following slightly strengthened form. The proof method is from an algorithm
of Dourisboure and Gavoille~\cite{dourisboure} to approximate tree-length. They used it to construct a tree-decomposition of a graph, with outer diameter $k$ say, where $\odw(G)\ge (k-1)/3$. We are going to extract something a little stronger from the same 
tree-decomposition: that $\glc(G)\ge k-1$ (which implies $\odw(G)\ge (k-1)/3$, by \ref{mainthm1}).

\begin{thm}\label{mainthm3}
{\bf Theorem:} Let $G$ be a graph; then $\glc(G)\ge 2 \lfloor \odw(G)/2\rfloor$. Moreover, if $G$ has a cycle,
then there is a cycle $C$ of $G$, and two
edge-disjoint paths $P,Q$ of $C$, both geodesic in $G$ and both with length $\lfloor \odw(G)/2\rfloor$, 
such that the loaded cycle $(C,E(P)\cup E(Q))$ is geodesic in $G$.
\end{thm}
\Proof If $\odw(G)\le 1$ the statement is clear, so
we assume that $\odw(G)\ge 2$.
Let $d\ge 2$ be an integer with $d\le \odw(G)$. We will prove that 
$\glc(G)\ge 2 \lfloor d/2\rfloor$ for all choices of $d$, which implies the theorem (even if $\odw(G)$ is infinite).  
Since the outer diameter-width of $G$ is the supremum of the 
outer diameter-width of its components, there is a component with outer diameter-width at least $d$, 
so we may assume that $G$ is connected. 
Choose $r\in V(G)$, and for $i\ge 0$ let $L_i$
be the set of all vertices $v$ such that $d_G(v,r)=i$. Thus the sets $L_0,L_1,\ldots$ are pairwise disjoint and have union $V(G)$.
For $i\ge 0$, say $u,v\in L_i$ are {\em equivalent} if there is a path between $u,v$ with interior in $L_i\cup L_{i+1}\cup\cdots$.
This is an equivalence relation; let $\mathcal{M}_i$ be the set of all equivalence classes. Let $\mathcal{M}$ be the union of the 
sets $\mathcal{M}_i$ over all $i\ge 0$. Thus, $\mathcal{M}$ is a collection of pairwise disjoint subsets of $V(G)$ with union $V(G)$.
Make a graph $T$ with vertex set $\mathcal{M}$, where distinct $A,B\in \mathcal{M}$ are adjacent if some vertex in $A$ is adjacent
in $G$ to some vertex of $B$. For each $A,B\in \mathcal{M}$, if $A,B$ are adjacent in $T$ then $A\in \mathcal{M}_i$
and $B\in \mathcal{M}_j$ for some $i,j\ge 0$ with $|j-i|=1$; and for each $i\ge 1$ and each $A\in \mathcal{M}_i$, there is a 
unique $B\in \mathcal{M}_{i-1}$ such that $A,B$ are adjacent in $T$, which we call the {\em parent} of $A$. We call $A$ a {\em child} of $B$.  Hence $T$ is a tree. 

For each $A\in \mathcal{M}$, let $W_A$ be the union of $A$ and its parent, if $A$ has a parent, and let $W_A:=A$ otherwise
(this only occurs when $A=\{r\}$).
\\
\\
{\bf (1) Claim: }{\em $(T,(W_A:A\in V(T)))$ is a tree-decomposition of $G$.}
\\
\\
Every edge of $G$ lies between $A$ and its parent, or has both ends in $A$, for some $A\in \mathcal{M}$, and in either case it has 
both ends in $W_A$; and therefore
the first two conditions in the 
definition of a tree-decomposition are satisfied. For the third, let $A,B,C\in V(T)$, where $B$ lies on the path of $T$ 
between $A,C$.
Let $v\in W_A\cap W_C$; we must show that $v\in W_B$. Choose $D\in \mathcal{M}$ with $v\in M$. Then the vertices $t\in V(T)$
with $v\in W_t$ are $D$ and its children in $T$; and so $A,C$ are both equal to or children of $D$. Consequently so is 
$B$, since it lies on the path of $T$ between $A,C$, and so $v\in W_B$. This proves (1).
\\
\\
{\bf (2) Claim:} {\em There exist $A\in \mathcal{M}$, and $u,v\in A$,
with $d_G(u,v)\ge d-1$.}
\\
\\
The outer diameter of $(T,(W_A:A\in V(T)))$ is at least $\odw(G)\ge d$, so
there exist
$A_0\in \mathcal{M}$, and $u_0,v_0\in W_{A_0}$, with $d_G(u_0,v_0)\ge d$. 
If $u_0,v_0\in A_0$ then the claim is true, so we may assume that $A_0$ has a parent $B$, and $u_0\in B$.
If $v_0\in B$ then again the claim is true, so we assume that $v_0\in A_0$. Thus $v_0$ has a neighbour $v\in B$, and 
$d_G(u_0,v)\ge d_G(u_0,v_0)-1\ge d-1$; and so the claim
is satisfied by $B,u_0,v$. This proves (2).

\bigskip

Let $A\in \mathcal{M}_i$, and 
$k=\lfloor d/2\rfloor\ge 1$.
Since $d_G(u,v)\ge d-1\ge 2k-1$, it follows that $i\ge k$. Let $P'$ be a shortest path of $G$ between $u,r$, and let $P$ be the subpath 
of $P'$ of length $k$ that contains $u$.
Let $p$ be the end of $P$ different from $u$ (so $p\in L_{i-k}$). Define $Q,q$ similarly, using $v$ in place of $u$.
Thus $P,Q$ are both geodesic in $G$.
Since $d_G(u,v)\ge 2k-1$, it follows that $P^*\cap Q^*= \emptyset$.
There is a path $R$
of $G$ between $u,v$ with interior in $L_i\cup L_{i+1}\cup\cdots$ (since $u,v$ are equivalent), 
and there is a path $S$ between $p,q$ with interior in $L_0\cupcup L_{i-k}$.
Since $i>i-k$, it follows that $R,S$ are vertex-disjoint.
Since $p$ is the only vertex of $P$ in $L_0\cup L_1\cupcup L_{i}$,
it follows that $V(P)\cap V(S)=\{p\}$, and similarly $V(P)\cap V(R)=\{u\}$, $V(Q)\cap V(S)=\{q\}$, and $V(Q)\cap V(R)=\{v\}$. Thus,
$P\cup Q\cup R\cup S$ is a cycle. Let $C$ be this cycle, and let $F:=E(P)\cup E(Q)$. 

Let 
$a,b\in V(C)$; we will show that $d_G(a,b)\ge d_{C,F}(a,b)$, and so $(C,F)$ is geodesic and the result holds.
Suppose that there is a path $T$ of $C$ between $a,b$ that is edge-disjoint from one of $P,Q$, and we assume it is edge-disjoint from $Q$ without loss of generality.
Let $a\in L_x$ and $b\in L_y$. Then
$|y-x|\ge |E(P)\cap E(T)|$; but $|y-x|\le d_G(a,b)$, and so
$$d_G(a,b)\ge |E(P)\cap E(T)|=d_{C,F}(a,b),$$
as required. Now suppose there is no such $T$, and hence one of $a,b$ belongs to $Q^*$, and one belongs to $P^*$;
so we assume that $a\in P^*$ and $b\in Q^*$.
Let $P_1,P_2$ be the subpaths of $P$ between $a$ and $u,p$ respectively, and define
$Q_1,Q_2$ similarly.
Let $\ell=|E(P_1)|+|E(Q_1)|$. Then 
$$\ell+d_G(a,b)\ge d_G(u,v)\ge 2k=\ell+|E(P_2)|+|E(Q_2)|\ge \ell+d_{C,F}(a,b)$$
and so again, $d_G(a,b)\ge d_{C,F}(a,b)$. 

Hence, $(C,E(P)\cup E(Q))$ is a geodesic loaded cycle,
and so 
$2\lfloor d/2\rfloor=2k\le \glc(G)$.
This proves \ref{mainthm3}.~\bbox

\section{Quasi-isometry}\label{sec:quasi}

We defined $(L,C)$-quasi-isometry in the first section.
If there exist $L,C$ such that $G$ admits an $(L,C)$-quasi-isometry to a tree, $G$ is a {\em quasi-tree}.
Every finite connected graph $G$ is  $(1,|G|)$-quasi-isometric to a tree, so all finite graphs are quasi-trees.
But infinite graphs may not be quasi-trees, and
(infinite) quasi-trees
are of substantial interest to
geometric group theorists~\cite{bala,bestvina}. It turns out that 
quasi-trees are the connected graphs with $\odw(G)$ finite.
We will show that
for a connected graph $G$, the following three statements are equivalent, and a bound in any one statement yields bounds for the other two:
\begin{itemize}
\item $\odw(G)$ is bounded;
\item there is an $(L,C)$-quasi-isometry to a tree with $L,C$ bounded;
\item there is an $(1,C)$-quasi-isometry to a tree with $C$ bounded.
\end{itemize}
(More exactly, if $\odw(G)$ is finite, then there is a $(1,6\odw(G))$-quasi-isometry to a tree; and if there is an  $(L,C)$-quasi-isometry to a tree, 
then $\odw(G)\le L(L+C+1)+C$.)
The equivalence of the second and third bullets here follows from Kerr's theorem \ref{kerr} that we mentioned earlier, 
but we will prove
it without assuming Kerr's theorem, because it seems to us that the equivalence is of sufficient interest to graph theorists that
it deserves a graph-theoretic proof.

Let us digress a little.
The equivalence of 
the second and third bullets above is striking, and one would naturally ask, how far does it extend? 
Is it confined to quasi-isometries to trees? The answer is no, but the question is awkward to make precise.
For it to have any content, we must be talking about 
quasi-isometries to graphs of some ``type'' (whatever that means!), not just to one graph; because for instance, 
a path $P$ of length $k$ is $(2,0)$-quasi-isometric to a path  $Q$ of length $2k$, but there is no $C$ (independent 
of $k$) such that there is a $(1,C)$-quasi-isometry from $P$ to $Q$. 
But at least the equivalence extends beyond trees: one can show that:
\begin{itemize}
\item for all $L,C$ there exists $C'$ such that 
if a finite graph $G$ is $(L,C)$-quasi-isometric to a cycle, then $G$ is $(1,C')$-quasi-isometric to a cycle (this is due to 
A. Georgakopoulos, in private communication);
\item for every integer $k\ge 1$, let $\mathcal{H}_k$ be the set of all finite connected graphs with no $K_{1,k}$ minor; then 
for all $L,C$ there exists $C'$ such that if a finite graph $G$ is $(L,C)$-quasi-isometric to a member of $\mathcal{H}_k$, then
$G$ is $(1,C')$-quasi-isometric to a member of $\mathcal{H}_k$ (this is an unpublished result of T. Nguyen, A. Scott and P. Seymour).
\end{itemize}
We do not know if this extends further. For instance, it seems to be open whether the equivalence holds for quasi-isometries to planar graphs.

Returning to the proof of the equivalence for trees, 
let us show first:
\begin{thm}\label{kerr2}
{\bf Theorem: }If $G$ is a connected graph with $\odw(G)$ finite, there is a $(1,6\odw(G))$-quasi-isometry from $G$ to a tree, and in particular
$\ad(G)\le 6\odw(G)$.
\end{thm}
\Proof
Let $k:=\odw(G)$, and let 
$(T,(B_t:t\in V(T)))$ be a tree-decomposition of $G$ with outer diameter $k$. We will show that  there is a 
$(1,6\odw(G))$-quasi-isometry from $G$ to a tree that is obtained from a subtree of $T$ by contracting and subdividing edges. Since $G$ is connected, we may assume
that $B_t\ne \emptyset$ for each $t\in V(T)$ (because the set of vertices $t\in V(T)$ with $B_t\ne \emptyset$ induces a subtree).
For all $s,t\in V(T)$, we denote the path of $T$ between $s,t$ by $T[s,t]$.

Choose $r\in V(T)$, and choose some vertex $\beta(r)\in B_r$.
For each $t\in V(T)$, let $\beta(t)$ be a vertex $v\in B_t$ with $d_G(v,\beta(r))$ minimum. 
For each edge $e=st\in E(T)$, let 
$$\ell(e):=|d_G(\beta(t),\beta(r))-d_G(\beta(s),\beta(r))|;$$ 
and
for each path $P$ of $T$, we
define $\ell(P)=\sum_{e\in E(P)}\ell(e)$.

For each edge $st$ of $T$, where $s$ is between $t$ and $r$, it follows that every path between $\beta(t)$ and $\beta(r)$
has a vertex in $B_s$, by \ref{intercept}, and so $d_G(\beta(s),\beta(r))\le d_G(\beta(t),\beta(r))$, because of the choice 
of $\beta(s)$.  Consequently
$\ell(T[t,r])=d_G(\beta(t),\beta(r))$ for each $t\in V(T)$.
For all $s,t\in V(T)$, we define $\ell(s,t)=\ell(T[s,t])$.

Now for each $v\in V(G)$, let $\phi(v)$ be some $t\in V(T)$ such that $v\in B_t$ (such a vertex exists from the definition
of a tree-decomposition). We will show that $d_G(u,v)$ and $\ell(\phi(u),\phi(v))$ differ by at most $6k$, for all $u,v\in V(G)$.
\\
\\
{\bf (1) Claim: }{\em If $v\in V(G)$, and $t\in V(T[\phi(v),r])$, then
\begin{align*}
\ell(\phi(v),t)&\le d_G(v,\beta(t)) +k, \text{ and}\\
d_G(v,\beta(t))&\le \ell(\phi(v),t) +3k.
\end{align*}
}
\noindent We have
$$\ell(\phi(v),t)= \ell(\phi(v),r)-\ell(t,r) =d_G(\beta(\phi(v)), \beta(r))-d_G(\beta(t),\beta(r))
\le d_G(\beta(\phi(v)), \beta(t))$$
by the triangle inequality. But $d_G(\beta(\phi(v)),v)\le k$ since $\beta(\phi(v))$ and $v$ both belong to $B_{\phi(v)}$, and so 
$$\ell(\phi(v),t)\le  d_G(\beta(\phi(v)), \beta(t))\le d_G(\beta(\phi(v)),v)+ d_G(v,\beta(t))\le d_G(v,\beta(t))+k.$$
This proves
the first statement.

For the second, let $P$ be a shortest path of $G$ between $\beta(\phi(v))$ and $\beta(r)$.
Then $P$ meets $B_t$, by \ref{intercept}; choose $x\in V(P)\cap B_t$.
Then
$$\ell(t,r)=d_G(\beta(t), \beta(r)) \le d_G(x, \beta(r))+k.$$
Since $x\in V(P)$, 
$$d_G(\beta(\phi(v)),x)= d_G(\beta(\phi(v)),\beta(r))-d_G(x,\beta(r))\le  \ell(\phi(v),r) - \ell(t,r)+k=\ell(\phi(v),t)+k.$$
Moreover,
$$d_G(v,\beta(t))\le d_G(v,\beta(\phi(v))) + d_G(\beta(\phi(v)), x)+d_G(x, \beta(t))\le  d_G(\beta(\phi(v)),x)+2k;$$ 
so $d_G(v,\beta(t))\le \ell(\phi(v),t) +3k$.
This proves (1).
\\
\\
{\bf (2) Claim:} {\em If $u,v\in V(G)$, then
\begin{align*}
\ell (\phi(u),\phi(v))&\le d_G(u,v) +4k, \text{ and}\\
d_G(u,v)&\le \ell(\phi(u),\phi(v)) +6k.
\end{align*}}
\noindent Let $t$ be the unique vertex of $T$ that belongs to all three of the paths that join two of $\phi(u),\phi(v),r$, and 
let $P$ be a shortest path of $G$ between $u,v$. It meets $B_t$ by \ref{intercept}, so
$$|E(P)|+2k\ge d_G(\beta(t),v) +d_G(\beta(t),u) \ge \ell(\phi(u),t) + \ell (\phi(v),t) -2k$$
by (1).
Hence $d_G(u,v)=|E(P)|\ge \ell(\phi(u),\phi(v))-4k$.
This proves the first statement.  For the second, by (1),
$$d_G(u,v)\le d_G(u,\beta(t))+d_G(v,\beta(t))\le \ell(t,\phi(u))+\ell(t,\phi(v)) +6k = \ell(\phi(u),\phi(v))+6k.$$
This proves (2).

\bigskip

Let $X:=\{\phi(v):v\in V(G)\}$, and let
$T'$ be the minimal subtree of $T$ with $X\subseteq V(T')$.
\\
\\
{\bf (3) Claim:} {\em For each $t\in V(T')$, there exists $v\in V(G)$ such that $\ell(\phi(v),t)\le 2k$.}
\\
\\
We may assume that $t\notin X$, and since $t\in V(T')$, $t$ belongs to a path of $T$ with both ends in $X$. Consequently
there are two components of $T\setminus \{t\}$ that both contain a vertex in $X$; and since $G$ is connected,
there are adjacent $u,v\in V(G)$ such that $\phi(u),\phi(v)$ belong to different components of $T\setminus \{t\}$. Hence
$t$ belongs to $T[\phi(u),\phi(v)]$. But $d_G(u,v)=1$, so by (2), $\ell(\phi(u),\phi(v))\le 4k+1$; and so
one of $\ell(\phi(u),t),\ell(\phi(v),t)$ is at most $2k$. This proves (3).

\bigskip

Now let $S$ be the tree obtained from $T'$ by contracting all edges $e$ with $\ell(e)=0$, and subdividing $\ell(e)-1$ times
(that is, replacing by a path of length $\ell(e)$) every edge $e$ with $\ell(e)>0$. If $t\in V(T')$, it has been identified
with other vertices of $T'$ under edge-contraction to form a vertex $\sigma(t)$ say of $S$. Thus for all $s,t\in V(T')$, $\ell(s,t)$
is the distance in $S$ between $\sigma(s),\sigma(t)$. From (2) and (3), it follows that
the map sending each $v\in V(G)$ to $\sigma(\phi(v))$ is a $(1,6k)$-quasi-isometry. 
Consequently $\ad(G)\le 6k$. This proves \ref{kerr2}.~\bbox

Next, we show:
\begin{thm}\label{quasitobags}
{\bf Theorem: }If there is an $(L,C)$-quasi-isometry from a graph $G$ to a tree,  then $G$ is connected and $\odw(G)\le L(L+C+1)+C$.
\end{thm}
\Proof
Let $\phi$ be an $(L,C)$-quasi-isometry from $G$ to a tree $T$. For all $u,v\in V(G)$, since $d_T(\phi(u),\phi(v))$ is finite, and 
$d_G(u,v)\le Ld_T(\phi(u),\phi(v))+C$, it follows that 
$d_G(u,v)$ is finite, and so $G$ is connected. For each $t\in V(T)$, let $B_t$ be the set of all $v\in V(G)$
such that $d_T(t,\phi(v))\le (L+C+1)/2$. Every vertex of $T$ within distance $(L+C+1)/2$ of both ends of a path $P$ of $T$
is within distance $(L+C+1)/2$ of every vertex of $P$; so if $t,t',t''\in V(T)$ and $t'$ belongs to the path
between $t,t''$ then $B_t\cap B_{t''}\subseteq B_{t'}$. Moreover, if $uv\in E(G)$, then $d_T(\phi(u),\phi(v))\le L+C$, and so
there exists $t\in V(T)$ within distance $\lceil (L+C)/2\rceil\le (L+C+1)/2$ from both $\phi(u),\phi(v)$, and hence $u,v\in V_t$.
It follows that $(T,(B_t:t\in V(T)))$ is a tree-decomposition. For each $t\in V(T)$, if $u,v\in B_t$ then
$$d_T(\phi(u),\phi(v))\le d_T(\phi(u),t)+d_T(\phi(v),t)\le L+C+1$$
and so $d_G(u,v)\le L(L+C+1)+C$. This proves \ref{quasitobags}.~\bbox

There were results already known that characterize when a graph is quasi-isometric to a tree.
The {\em bottleneck constant} of a graph $G$ is the least integer $\Delta$ such that 
if $P$ is a geodesic path of $G$ between $u,v$, of even length and with middle vertex $w$, then every path between $u,v$ contains a vertex that has 
distance at most $\Delta$ from $w$.
A theorem of Manning for geodesic metric spaces implies: 
\begin{thm}\label{manning2}
{\bf Theorem \cite{manning}: }For all $L\ge 1$ and $C\ge 0$, there exists $\Delta$ such that, for all graphs $G$,
 if there is an $(L,C)$-quasi-isometry from $G$ to a tree, then
$G$ has bottleneck constant at most $\Delta$. Conversely, for all $\Delta$ there exist $L\ge 1$ and $C\ge 0$ such that, 
for all graphs $G$, if $G$ has bottleneck constant at most $\Delta$, then there is an $(L,C)$-quasi-isometry from $G$ to a tree.
\end{thm}
We observe also:
\begin{thm}\label{bottle}
{\bf Lemma: }If $G$ has bottleneck constant $\Delta$, then $\glc(G)\ge 2\Delta$.
\end{thm}
\Proof
Suppose that $G$ has bottleneck constant $\Delta\ge 1$. The minimality of $\Delta$ implies that
there exist $u,v,w\in V(G)$ with $d_G(u,w)=d_G(w,v)=d_G(u,v)/2$, and there is a path
$S$ between $u,v$ such that all its vertices have distance at least $\Delta$ from $w$. Let $P_0$ be a path between $u,w$
of length $d_G(u,w)$, and define $Q_0$ similarly with ends $v,w$. 
Since $u\in V(S)$, it follows that $P_0$ has length at least $\Delta$; choose $p\in V(P_0)$ such that the subpath $P$ 
of $P_0$ between $p,w$
has length $\Delta$. Choose $q\in V(Q)$ and $Q$ similarly. Then the union of the path of $P_0$ between $p,u$, the path $S$, and 
the path of $Q_0$ between $v,q$, contains a path $R$ between $p,q$ such that all its vertices have distance at least $\Delta$ 
from $w$. If $x\in V(P)$ and $y\in V(R)$, then $d_G(x,y)+d_G(x,w)\ge \Delta$; but $d_G(x,w)+d_G(p,x) = \Delta$, and so 
$d_G(x,y)\ge d_G(p,x)$. A similar statement holds if $x\in V(Q)$; and it follows that 
 $P\cup Q\cup R$ is a cycle, and $(P\cup Q\cup R,E(P\cup Q))$ is a geodesic loaded cycle, with load $2\Delta$. This proves \ref{bottle}.~\bbox

So one can deduce from Manning's theorem and \ref{bottle}
that if $\glc(G)$ is bounded, then there is an $(L,C)$-quasi-isometry from $G$ to a tree with $L,C$ bounded, and hence,
from \ref{quasitobags}, that $\odw(G)$ is bounded. Our result \ref{mainthm} says the same, but with more
explicit control over the bounds.
We can also deduce a version of Manning's theorem (for graphs: Manning's theorem is really for metric spaces) from our result:
\begin{thm}\label{tomanning}
{\bf Theorem: }If a graph $G$ has bottleneck constant $\Delta$, then $\odw(G)\le 4\Delta+3$, and hence there is a $(1,24\Delta+18)$-quasi-isometry to a tree.
\end{thm}
\Proof
Suppose that there is a 
geodesic loaded cycle
$(C,F)$ with load $\ge 2\Delta+2$, such that $F$ is the edge-set of a geodesic path $P$. We may assume that $P$ has length 
$2\Delta+2$.
Let $w$ be its middle vertex. Then $C\setminus P^*$ is a path 
between the ends of $P$, and for each of its vertices $v$, $d_G(v,w)\ge d_{C,F}(v,w)= \Delta+1$, contrary to the definition 
of bottleneck constant. Thus there is no such $(C,F)$.

Hence, from the final statement of \ref{mainthm3}, it follows that $\lfloor \odw(G)/2\rfloor \le 2\Delta+1$,
and so $\odw(G)\le 4\Delta+3$. Applying \ref{kerr2} now proves \ref{tomanning}.~\bbox

Finally, there are further, similar, characterizations in~\cite{chepoi, george}.

\section{McCarty's conjecture}\label{sec:mccarty}

Rose McCarty (private communication) suggested a different condition that might characterize when there is a tree-decomposition
of small outer diameter, as follows. Let us say a graph has {\em McCarty-width} $k$ if $k\ge 0$ is minimum such that the following holds:
for every three vertices $u,v,w$ of $G$, there is a vertex $x$, such that if $X$ denotes the set of all vertices that have distance
at most $k$ from $x$, then no component of $G\setminus X$ contains two of $u,v,w$. Let $\mcw(G)$ denote the McCarty-width of $G$.
McCarty suggested that $\odw(G)$ is small 
if and only if $\mcw(G)$ is small. This turns out to be true, because of the following:

\begin{thm}\label{mccarty}
{\bf Theorem: }Let $G$ be a graph. Then $(\odw(G)-3)/6\le \mcw(G)\le \odw(G)$.
\end{thm}
\Proof
We show first that $\mcw(G)\le \odw(G)$.
Let $(T,(B_t:t\in V(T)))$ be a tree-decomposition of $G$ with outer diameter equal to $\odw(G)$. Now let $u,v,w\in V(G)$,
and choose $t_1,t_2,t_3\in V(T)$ with $u\in B_{t_1}$, $v\in B_{t_2}$ and $w\in B_{t_3}$. Let $t$ be the unique vertex of $T$ that
belongs to each of the three paths of $T$ that join two of $t_1,t_2,t_3$. Let $x\in B_t$; then every path of $G$ between two 
of $u,v,w$ contains a vertex of $B_t$, by \ref{intercept}, and all such vertices have distance at most $\odw(G)$ from $x$. Hence $\mcw(G)\le \odw(G)$.

For the other inequality, suppose that $\glc(G)\ge  6\mcw(G)+3$, and 
choose a geodesic loaded cycle $(C,F)$ of $G$ with $|F|=\glc(G)$.
Choose three distinct vertices $u,v,w\in V(C)$, such that each of $d_{C,F}(u,v), d_{C,F}(u,w),d_{C,F}(v,w)$ is at least
$2\mcw(G)+1$. Let $C_{u,v}$ be the path of $C$ between $u,v$ not containing $w$, and define
$C_{u,w},C_{v,w}$ similarly. From the definition of McCarty-width, there is a vertex $x$, 
such that if $X$ denotes the set of all vertices that have distance
at most $\mcw(G)$ from $x$, then no component of $G\setminus X$ contains two of $u,v,w$. Hence some vertex of $C_{u,v}$
belongs to $X$, say $c_{u,v}$, and define $c_{u,w}, c_{v,w}$ similarly. Since $c_{u,v}, c_{u,w}$ both have distance
at most $\mcw(G)$ from $x$, they have distance at most $2\mcw(G)$ from each other, and so $d_{C,F}(c_{u,v}, c_{u,w})\le 2\mcw(G)$.
Since $C_{v,w}$ contains at least $2\mcw(G)+1$ edges in $F$, it follows that
the path ($P_u$ say) of $C$ between $c_{u,v}, c_{u,w}$ that does not include $c_{v,w}$ contains at most $2\mcw(G)$ edges
of $F$.
The same holds for the other two pairs of $c_{u,v}, c_{u,w}, c_{v,w}$; define $P_v,P_w$ similarly. 
But every edge in $F$ belongs to one of $P_u,P_v,P_u$, and so $|F|\le 6\mcw(G)$, a contradiction. This proves that
$\glc(G)\le  6\mcw(G)+2$, and since $\glc(G)\ge \odw(G)-1$ by \ref{mainthm},
it follows that 
$(\odw(G)-3)/6\le \mcw(G)$. This proves \ref{mccarty}.~\bbox
\section{On spanning tree distortion}\label{sec:counterex}

There was another candidate that we hoped would characterize when $\odw(G)$ was small,
as follows. We know that if $\odw(G)$ is small, there is a $(1,C)$-quasi-isometry $\phi$ to a tree $T$, 
and we might hope that $T$ can be chosen to be a spanning tree of $G$, and $\phi$ the identity function.
Let us say the {\em cycle-distortion} of a spanning tree $T$ of $G$ is the maximum, over all edges $uv$ of $G$,
of the length of the path of $T$ between $u,v$. If $G$ admits a spanning tree with cycle-distortion $d$, then $\idw(G)\le 2d$
(use the same tree, with the bag for vertex $t$ a ball of $T$ with radius $d$ around $t$), so one might 
hope for a converse, to give a characterization, at least for connected graphs $G$. But this is not the case,
because of the following (a closely-related result appears in section 6 of~\cite{kratsch}):

\begin{thm}\label{counterex}
{\bf Theorem: }There is a connected graph $G$ with $\idw(G)=1$, such that every spanning tree has large cycle-distortion.
\end{thm}
\Proof Let $D_1$ be a cycle of length three, drawn in the plane: so its outer boundary is (trivially) a cycle $C_1=D_1$. 
Inductively, for
$i\ge 2$, we assume that $D_{i-1}$ is drawn in the plane with its outer boundary a cycle $C_{i-1}$: let $D_i$ be obtained from $D_{i-1}$ by adding
a new vertex $z_{uv}$ for each edge $uv$ of $C_{i-1}$, adjacent to $u$ and to $v$, drawn outside $C_{i-1}$ such that
the outer boundary $C_i$ of $D_i$ is formed by these new edges. (Thus, each $D_i$ is a finite subgraph of the ``Farey graph''.)
Let $k$ be a large integer.
We claim that $\idw(D_k)=1$, and every spanning tree of $D_k$ has cycle-distortion at least $k+1$.

To see the first claim, note that $D_k$ is a chordal graph, and therefore admits a tree-decomposition $(T,(B_t:t\in V(T)))$,
where each $B_t$ is a clique of $D_k$, and so has inner diameter-width $1$. For the second claim, 
let $T$ be a spanning tree of $D_k$. Every vertex of $D_k$ belongs to $C_k$; and for every edge $e=uv$ of $D_k$, not an edge of 
$C_k$, we observe that $\{u,v\}$ separates $D_k$ into exactly two components. A {\em triangle} is a cycle of length three. Every
triangle of $D_k$ is the boundary of a region of the drawing, and every finite region (that is, every region except the 
infinite region)  has boundary a triangle. If
$e=uv$ is an edge of a triangle $\Delta$, we say that $e$ is {\em $\Delta$-bad} if the $u\DD v$ path of $T$ is
vertex-disjoint from the component of $D_k \setminus \{u,v\}$ that contains the third vertex of $\Delta$. (In
particular, if $e \in E(T)$ then $e$ is $\Delta$-bad.) If all three edges of a triangle $\Delta$
are $\Delta$-bad, then the union of the corresponding three paths is a cycle of $T$,
which is impossible. Let $\Delta_1:=C_1$. At least one edge $e_1=u_1v_1$ of $\Delta_1$ is not $\Delta_1$-bad; let $P_1$ be
the path of $T$ between $u_1,v_1$. Thus $P_1$ contains all three vertices of $\Delta_1$, and has length at least two. 
Let $\Delta_2$ be the other triangle
containing $e_1$. Thus $e_1$ is $\Delta_2$-bad, and so some other edge $e_2=u_2v_2$ of
$\Delta_2$ is not $\Delta_2$-bad. The $u_2\DD v_2$-path $P_2$ of $T$ contains both ends of $P_1$ and so contains $P_1$,
and hence has length at least three.
Repeating the argument
inductively, we obtain a nested sequence $P_1, P_2,\ldots P_k$ of paths of $T$, each with
adjacent ends, where each $P_i$ has length at least $i+1$; and the cycle-distortion of $T$ is at least the length of all these paths. This proves \ref{counterex}.~\bbox

\section*{Acknowledgements}
We would like to thank Alex Scott and Rose McCarty, who worked with us on an earlier (false) conjectured solution to this problem. 
We would also like to thank Alex for his help with the manuscript, and Rose for the suggestion of McCarty-width.
And we are grateful to the referees for their careful, helpful reports; to one referee in particular,
who introduced us to quasi-isometry; to Carla Groenland, who told us about
tree-length; and to Agelos Georgakopoulos, for help in understanding quasi-isometry better.

\end{document}